\documentclass{report}
\usepackage{amssymb}
\usepackage{amsthm}
\usepackage{eufrak}
\usepackage[all]{xy}
\usepackage{makeidx}

\makeindex
\title{Notes on Cohomology}
\author{Luis Arenas-Carmona.
\footnote{Supported by Fondecyt, proyecto No. 3010018,
and the Chilean Catedra Presidencial in Number Theory.}\\
\\
\small\textit{Universidad de Chile,}\\
\small\textit{Facultad de Ciencias.}\\
\small{Casilla 653, Santiago,Chile.}\\
\small{learenas@uchile.cl}}
\date{}

\theoremstyle{plain}
\newtheorem{prop}{Proposition}[section]

\newtheorem{lem}[prop]{Lemma}
\newtheorem{cor}{Corollary}[prop]
\theoremstyle{definition}
\newtheorem{ex}[prop]{Example}
\newtheorem{dfn}[prop]{Definition}
\newtheorem{rmk}[prop]{Remark}

\newcommand\stab{\textnormal{Stab}}
\newcommand\End{\textnormal{End}}
\newcommand\Aut{\textnormal{Aut}}
\newcommand\Hom{\textnormal{Hom}}
\newcommand\smallgen{\textnormal{\tiny{gen}}}

\newcommand\im{\textnormal{im}}

\newcommand\ckgen[3]{C_{\smallgen}(#1,#2,#3)}

\newcommand\alge{\mathfrak{A}}

\newcommand\ad{\mathbb A}

\newcommand\Stab{\textnormal{Stab}}

\newcommand\spaan{\textnormal{span}\ }

\newcommand\la{\Lambda}
\newcommand\bark{\bar{k}}

\newcommand\uno{\{1\}}

\newcommand\oink{\mathcal O}

\newcommand\Q{\mathbb Q}

\newcommand\reali{\mathbb R}
\newcommand\compleji{\mathbb C}
\newcommand\enteri{\mathbb Z}

\newcommand\gal{\mathcal G}
\newcommand\hal{\mathcal H}
\newcommand\galbark{{\mathcal G_{\bark/k}}}

\newcommand\galbarkk{\galbark}
\newcommand\galkk{{\mathcal G_{K/k}}}
\newcommand\galkvkv{{\mathcal G_{K_w/k_v}}}
\newcommand\galv{{\mathcal G_w}}
\newcommand\tn{\mathfrak{T}}
\newcommand\dete{\textnormal{det}}
\newcommand\coker{\textnormal{coker}}

\newcommand\vaa{\longrightarrow}

\newcommand\grupi{\mathbb G}

\newcommand\ideala{\mathcal A}
\newcommand\idealp{\wp}
\newcommand\idealP{\mathcal P}

\newcommand\orth[3]{\oink_{#1,#2}(#3)}

\newcommand\orthi[4]{\oink_{#1,#2}^{#4}(#3)}

\newcommand\orthomega{\oink_n(Q)}

\newcommand\orthkbar{\orth n{\bark}Q}

\newcommand\linei[4]{GL_{#2}^{#4}(#3)}
\newcommand\linear[3]{\linei {#1}{#2}{#3}{}}
\newcommand\linek{\linear nkV}
\newcommand\lineE{\linear nEV}
\newcommand\lineK{\linear nKV}

\newcommand\lineklam{\linei nkV\la}
\newcommand\lineKlam{\linei nKV\la}

\newcommand\lineKvlam{\linei n{K_w}V\la}
\newcommand\lineomega{GL(V)}
\newcommand\linekbarlam{\linei n{\bark}V\la}

\newcommand\speciali[3]{SL_{#1}^{#3}(#2)}
\newcommand\specia[2]{\speciali {#1}{#2}{}}

\newcommand\specialK{\specia KV}

\newcommand\specialKlam{\speciali KV\la}

\newcommand\specialomega{SL(V)}

\newcommand\tensor{\bigotimes_{\oink_k}}

\newcommand\lk{\mathcal{L}_{\textnormal{\tiny def}}}
\newcommand\lfree{\mathcal{L}_{\textnormal{\tiny fr}}}
\newcommand\lvk{\mathcal{L}^V}
\newcommand\lkvk{\mathcal{L}_{\textnormal{\tiny def}}^V}
\newcommand\lfreevk{\mathcal{L}_{\textnormal{\tiny fr}}^V}

\newcommand\Lk[3]{\mathcal{L}_{\textnormal{\tiny def}}(#1,#2,#3)}
\newcommand\Lfree[3]{\mathcal{L}_{\textnormal{\tiny fr}}(#1,#2,#3)}
\newcommand\Lvk[3]{\mathcal{L}^V(#1,#2,#3)}
\newcommand\Lkvk[3]{\mathcal{L}_{\textnormal{\tiny def}}^V(#1,#2,#3)}
\newcommand\Lfreevk[3]{\mathcal{L}_{\textnormal{\tiny fr}}^V
(#1,#2,#3)}

\newcommand\cuaternioni{\mathbb H}
\newcommand\matrici{\mathbb M}
\newcommand\ksobrekn{k^*/(k^*)^n}

\newcommand\gtilde{\widetilde{G}}

\newcommand\gggr[2]{G_{#1}^{#2}}

\newcommand\gkla{\gggr k\la}
\newcommand\gakla{\gggr {\ad_k}\la}
\newcommand\gKla{\gggr K\la}

\newcommand\gKlav{\gggr {K_w}\la}

\newcommand\gK{G_K}
\newcommand\gk{G_k}

\newcommand\gkv{G_{k_v}}
\newcommand\gak{G_{\ad_k}}

\begin{document}
\maketitle



\chapter{Introduction.}

Galois cohomology is a fundamental tool for the
classification of certain algebraic structures.
To be precise, let $k$ be a field, $G$ a linear algebraic group
acting on a space $V$, both defined over $k$.
It is known \cite{kneser2}, that if $G$ is
 defined as
 the set of automorphisms of a tensor $\tau$ on $V$,
 e.g., a quadratic form or an algebra structure,
 the cohomology set $H^1(K/k,G_K)$ classifies the
 $K/k$-forms of $\tau$,
 i.e., those tensors of the same
 type also defined over $k$ that
 become isomorphic to $\tau$
over the larger field $K$ (\S \ref{kkforms}).
Results of this type, however, hold in much more general setings.
In this notes, we give the general facts about cohomology that allow
the use of cohomology sets for classification, and give examples
of aplications to many parts of field theory and number theory.
In particular, we devote a whole chapter to the study of the relation between
lattices and cohomology.

Such a theory is already hinted at in \cite{Rohlfs}.  In this
reference, two finiteness results are proven. The first one deals with
the finiteness of the local cohomology set $H^1(\gal_w,\Gamma_w)$,
for an arithmetically defined group $\Gamma$.
Notations are as in \cite{Rohlfs}.
The second one deals with the finiteness of the kernel of the map
\[
H^1(\gal,\Gamma)\vaa \prod_{v\textnormal{ place of
}k}H(\gal_{w(v)},\Gamma_{w(v)}),
\]
where we have fixed a place $w(v)$
of $K$ dividing each place $v$ of $k$.
It is the proof of the second result which requires
expressing the given kernel in terms of the set of double cosets
\[
G_k\backslash G_{\ad_k}/\prod_w\Gamma_w
\]
 (see corollary 3.3 in \cite{Rohlfs}).
These double cosets  are the same ones that
classify the classes of lattices in a genus. This relation is pursued
in chapter (crossreference).

\subsection{Notations}\label{notations}

In all of this notes, $k,K,E$ denote number or local fields of
characteristic $0$, or algebraic extensions of them.
If $k$ is a number field, $\Pi(k)$ denotes the set of
places of $k$.
\begin{rmk}\label{omegagroups}
By an algebraic group, we mean a linear algebraic
group. All algebraic groups are assumed to be subgroups
of the general linear group of a vector space $V$,
of finite dimension, over
a sufficiently large
algebraically closed
field $\Omega$ of characteristic $0$.
We assume that all localizations of number fields
inject into $\Omega$.
$G$ denotes an algebraic group over $\Omega$.
 $\lineomega,\specialomega$ denote the general and
special linear groups over $\Omega$.
When we work over a
fixed local or number field $k$,
we say that $G$ is defined over $k$ if the equations defining $G$
have coefficients in $k$ (see section 2.1.1 in \cite{Pla}).
 \emph{This is the case for all groups considered
here}.
For any field $E$, $k\subseteq E\subseteq\Omega$, we write $G_E$ for the set
of $E$-points of $G$, e.g., if $G=\lineomega$, the set of $E$ points
is denoted $\lineE$. The same conventions apply to spaces and algebras.
All spaces and algebras are assumed to be finite dimensional.

Exceptions to this rule are the multiplicative and additive groups.
We denote $\grupi_m=\Omega^*$, $\grupi_a=\Omega$ when considered as
algebraic groups. For the set of $k$-points
we write $k^*$,$k$. Instead of $(\grupi_m)_k,(\grupi_a)_k$.
\end{rmk}

The orthogonal group of a quadratic form $Q$ on $V$ is written
$\oink_{n}(Q)$ or $\oink_{n}(Q,V)$, where $n=\dim_\Omega(V)$.
The set of $E$-points is denoted $\oink_{n,E}(Q)$.

The field on which a particular lattice is defined is always written as
a subindex. If $K/k$ an extension of local or number fields
and $\la_k$ is a lattice in $V_k$, $\la_K$ denotes the $\oink_K$-lattice
in $V_K$ generated by $\la_k$.

If $G$ is an algebraic group  acting on a space $V$, both
defined over $k$,
and $\la_k$ is a $\oink_k$-lattice on $V_k$, the stabilizer of $\la_k$
in $G_k$ is denoted $\gkla$. If $G=\lineomega$, this set is denoted
$\lineklam$. Similar conventions apply to
special linear or orthogonal groups.

\begin{rmk}\label{fixw}
Whenever $K/k$ is a Galois extension of a number field $k$, and
 $v$ a place of $k$,
$w$ denotes a place of $K$ dividing $v$. We assume that one fixed such $w$
 has been chosen for every $v$. This convention is also applied for
infinite extension, e.g., $K=\bark$.
\end{rmk}

\begin{rmk}\label{galoisgroups}
$\galkk$ denotes the Galois group of the extension $K/k$.
If there is no risk of confusion, we write simply $\gal$.
If $K$ is not specified, we assume $K=\bar{k}$.
If $k$ is a number field and $v\in\Pi(k)$,
we also use the notation $\galv=\galkvkv$.
\end{rmk}

If $\Gamma$ is a group acting on a set $S$, $S/\Gamma$ denotes the set of orbits
and $S^\Gamma$ the set of invariant points. The action of $\gamma\in\Gamma$
is denoted $s\mapsto s^\gamma$, for $s\in S$.

\chapter{Cohomology and classification}

In this chapter we introduced the basic results that are required to connect cohomology and classification.
The results in this section are found in chapter 1 in \cite{kneser},
and p. 13-26 in \cite{Pla}.
\index{cohomology}

\begin{dfn}
Let $\gal$ be a finite group, and let $A$ a group provided with a $\gal$-action.
 $H^1(\gal,A)$ is defined as the quotient
\[
H^1(\gal,A)=\{\alpha:\gal\mapsto A|\alpha(hg)=\alpha(h)
\alpha(g)^h\}/\equiv,
\]
where $\alpha\equiv\beta$ if and only if there exists
$a\in A$  such that $\alpha(g)=a^{-1}\beta(g)a^g$ for all $g\in\gal$.
If $\gal$ acts trivially on $A$, then $H^1(\gal,A)\cong\Hom(\gal,A)/A$,
where $A$ acts on $\Hom(\gal,A)$ by conjugation.
In what follows we write $\alpha_g$ instead of $\alpha(g)$.
\end{dfn}

In case that $A\subseteq B$ is a subgroup, there is a long exact sequence
\[
 0\vaa A^\gal \vaa B^\gal \vaa (B/A)^\gal \vaa H^1(\gal,A) \vaa
 H^1(\gal,B),
\]
and furthermore, under the natural action
\footnote{ This result is not found in \cite{kneser}, but can be
found in \cite{Pla} p.22.}
of $B^\gal$ on $(B/A)^\gal$,
\begin{equation}\label{clasificador}
(B/A)^\gal/B^\gal \cong \ker(H^1(\gal,A) \vaa H^1(\gal,B)).
\end{equation}

To simplify notations, in all that follows we assume that
whenever a sequence of pointed sets
\[
\dots\vaa U\vaa V\vaa W\vaa X\vaa Y \vaa Z
\]
is written, $X,Y,Z$ denote pointed sets,
$W,V,U,\dots$ denote groups,
and $W$ acts on $X$ with
\[
X/W\cong\ker(Y\vaa Z).
\]

If $A$ is normal in $B$, we have in the sense just described
\begin{equation}\label{iso2}
\xymatrix{ 0\ar[r] & A^\gal\ar[r] & B^\gal\ar[r]  &
(B/A)^\gal \ar`r[d]`[lll]`[llld]`[llldr] [dll] \\
{}& H^1(\gal,A)\ar[r] & H^1(\gal,B)\ar[r] & H^1(\gal,B/A). }
\end{equation}
%

In case $A$ is central in $B$,
 the higher
order cohomology groups for $A$ are also defined, and we have a long exact
 sequence
\[
\xymatrix{  0\ar[r]& A^\gal\ar[r]& B^\gal\ar[r]&
(B/A)^\gal\ar`r[d]`[lll]`[llld]`[llldr] [dll]  \\
{} & H^1(\gal,A)\ar[r]& H^1(\gal,B)\ar[r]&
H^1(\gal,B/A)\ar`r[d]`[lll]`[llld]`[llldr] [dll] \\
{}& H^2(\gal,A). &{}&{}}
\]

Finally, if $A$ and $B$ are both Abelian this sequence extends
 to cohomology
of all orders \cite{Serre}, \cite{kenneth}.
All results of this section can be extended via direct limits
to profinite groups
acting continuously on discrete groups
(\cite{Serre}, p. 9 and 42).

\section{The general classification principle}\label{ReviewDefinitions}

A $(\gal,G)$-space  is a set $X$ provided with  both,
a $\gal$-action and a $G$-action (denoted $*$) satisfying
$$g^\sigma\ast x^\sigma=(g\ast x)^\sigma$$ for
$x\in P$, $g\in G$, and $\sigma\in\gal$. The space $X$ is said to be free, transitive, etc,
if the corresponding $G$-action has any of these properties.
The most important aplication for us of the results in the preceeding
 section is the following:

\begin{prop}\label{classificationlaw}
Let $X$ be a transitive $(\gal,G)$-space
for all $x\in X$, $g\in G$, and $\sigma\in\gal$.
Let $x_0\in X$, and let $H=\Stab_G(x_0)$.
Then, $X^\gal$
is in one-to-one correspondence with the elements of the cohomology
set $\ker\big(H^1(\gal,H)\rightarrow H^1(\gal,G)\big)$.
\end{prop}

\subparagraph{Proof.}
Since $X$ is isomorphic to $G/H$ as $\gal$-modules, the result follows from
(\ref{clasificador}).\qed

A principal homogeneous $G$-space is a  $(G,\gal)$-space where the $G$ action is transitive and free.
Given any principal homogeneous $G$-space $P$ and any base elements
$x\in P$, there exists a unique $G$-map (but not a $\gal$-map, unless $x$ is invariant)
 $\phi:G\rightarrow P$ that sends $g$ to $gx$. In particular,
we can identify $P$ with $G$ as a set.
Classifying principal homogeneous $G$-spaces is, therefore, equivalent to
classifying twisted actions on $\gal$ on $G$ that turn the set $G$ into a principal hemogeneous
$G$-space.

Let $X$ be the group of all maps from $\gal$ to $G$. Then $X$ is a $\gal$-group with
an action $f\mapsto f^\sigma$ satisfying $f^\sigma(\lambda)=f(\lambda\sigma^{-1})$.
A twisted semiaction of $\gal$ on $G$ is a map $\rho:G\times\gal\rightarrow G$
which satisfy $\rho(hg,\sigma)=h^\sigma\rho(g,\sigma)$ and $\rho(1,1)=1$.
A twisted semiaction $\rho$ is a twisted action if it satisfies the relation
$$\rho(\rho(g,\sigma),\lambda)=\rho(g,\sigma\lambda),$$
for all $g$ in $G$ and all $\sigma$ and $\lambda$ in $\gal$.
The group $X$ acts transitively on the set $T$ of twisted semiactions by
$$(\alpha\cdot\rho)(g,\sigma)=\rho(g\alpha(\sigma)^{-1},\sigma)\alpha(1),$$
for all $\alpha$ in $X$, all $\rho$ in $T$, all $g$ in $G$, and all $\sigma$ in $\gal$.
The stabilizer of $\rho_0$, where $\rho_0(g,\sigma)=g^\sigma$,
is the set of maps satisfying $\alpha(\sigma)^{\sigma}=\alpha(1)$.
This maps form a subgroup $G'$ of $X$ which is isomorphic to $G$ as a $\gal$-group.
The quotient set $X/G'\cong T$ has a natural $\gal$ action which is traslated to
to an action on $T$ as $\rho^{\lambda^{-1}}(1,\sigma)=\rho(1,\sigma\lambda)^{\lambda^{-1}}
\rho(1,\lambda)^{-\lambda^{-1}}$. With this action, a twisted semiaction $\rho$ is invariant if and only if
$$\rho(1,\sigma\lambda)=\rho(1,\sigma)^\lambda\rho(1,\lambda)=\rho(\rho(1,\sigma),\lambda),$$
and premultiplying both sides by $g^{\sigma\lambda}$ we see that a $\gal$-invariant action is the same as
a twisted action. The following lemma follows easily from
proposition \ref{classificationlaw}.

\begin{lem}\label{last}
The set of twisted actions on $G$ is in natural correspondence, up to isomorphisms of $(G,\gal)$-actions,
 with the kernel of the map $H^1(\gal,G)\rightarrow H^1(\gal,X)$.\qed
\end{lem}

\begin{prop}[Shapiro's Lemma]
Let $\hal$ be a subgroup of $\gal$.
Let $G$ be a $\hal$-group and let $G'$ be the set of maps $\phi:\gal\rightarrow G$
such that $\phi(\sigma\lambda)=\phi(\sigma)^\lambda$ for all $\lambda\in \hal$.
Then $H^1(\gal,G')\cong H^1(\hal,G)$.
\end{prop}

\begin{cor}\label{corl}
In the notations of lemma \ref{last} $H^1(\gal,X)=\uno$.\qed
\end{cor}

The following result follows from lemma \ref{last} and corollary \ref{corl}

\begin{prop}[(\cite{Serre}, p. 44)]
The set of principal homogeneous $G$-spaces up to isomorphism
is in one-to-one correspondence with the elements of the cohomology
set  $H^1(\gal,G)$.\qed
\end{prop}

In most applications of the results in this chapter $\gal$ is
the Galois group $\galkk$
of a  possibly infinite Galois extension $K/k$, where $k$ is a
local or number field.
The subgroups $A,B,\dots$
are usually groups of algebraic or arithmetical nature.

\chapter{The generalized Hilbert theorem 90}

\section{Rings, units and cohomology}

In this section, $A$ is a ring provided with a $\gal$ action. A left $(A,\gal)$-module
is an $A$-module $B$ provided with a $\gal$-action, and satisfying
$(ab)^\sigma=a^\sigma b^\sigma$ for all $a\in A$, $b\in B$, and $\sigma\in\gal$.

An element $b\in B$ is a generator if $Ab=B$. It is a regular generator if
$ab=0$ impplies $a=0$. In particular, if $b\in B$ is a regular generator,
the map $a\mapsto ab$ is an isomorphism of $A$-modules between $A$ and $B$.
Assume that $b$ is a regular generator of $B$. Then
an element $ab\in B$ is a regular generator if and only if $a\in A^*$.
We say that the $(A,\gal)$-module is principal if it has a regular generator.

\begin{prop}
The cohomology set $H^1(\gal,A^*)$ classifies the set of principal left
$(A,\gal)$-modules. The distinguished point of  $H^1(\gal,A^*)$
corresponds to the modules that have an invariant regular generator.
\end{prop}

\subparagraph{Proof.}
Let $X$ be the group of all maps from $\gal$ to $A^*$.
As in \S \ref{ReviewDefinitions}, we define a twisted semiaction as a map
$\rho:\gal\times A\rightarrow A$ satisfying:

a) $\rho(aa',\sigma)=a^\sigma\rho(a',\sigma)$.

b) $\rho(1,\sigma)$ is a unit for all $\sigma\in\gal$.

A twisted action as a twisted semiaction that is an action.
Any principal left $(A,\gal)$-module is isomorphic to $A$ with a twisted
action. The group $X$ acts transitively on the set of twisted semiaction
by $(\alpha\cdot\rho)(a,\sigma)=\rho(a\alpha(\sigma),\sigma)\alpha(1)^{-1}$,
and the stabilizer of $\rho_0$, where $\rho_0(a,\sigma)=a^\sigma$, is the
set of maps satisfying $\alpha(\sigma)^\sigma=\alpha(1)$ for all $\sigma\in\gal$.
Now the proof follows as in \S \ref{ReviewDefinitions}.
\qed

\begin{ex}
Let $A=\enteri$ with the trivial $\gal$-action. Then $H^1(\gal,\enteri^*)
=\Hom(\gal,\{1,-1\})$, hence there exists a free left $(\gal,\enteri)$-module
with no invariant generator if and only if $\gal$ has a normal subgroup of index $2$.
The reader can easily check this result independently.
\end{ex}

\begin{ex}
More generally, if $A$ has the trivial $\gal$-action. Then $H^1(\gal,\enteri^*)
=\Hom(\gal,A^*)/\cong$, where $\cong$ denotes the conjugacy relation.
The module associated to the homomorphism $\alpha:\gal\rightarrow A^*$
has the action $(ab)^\sigma=a\alpha(\sigma)b$.
\end{ex}

\begin{ex}
Let $K/k$ be a finite extension of local or number fields with
Galois group $\gal$. Let $\oink_K$ and $\oink_k$ be the respective
rings of integers. Then the exact sequence
$\oink_K^*\hookrightarrow K^*\twoheadrightarrow P_K$, where $P_K$
is the group of principal fractional ideals of $K$, shows that the
kernel of $H^1(\gal,\oink_K^*)\rightarrow H^1(\gal,K^*)$
corresponds to invariant ideals in $P_K$ modulo ideals in $P_k$
(up to a suitable identification). Since we prove later in this
chapter that $H^1(\gal,K^*)=1$, then all free left
$(\oink_K,\gal)$-module of rank $1$ are of this form. If the
extension $K/k$ is unramified, one can prove that $P_K^\gal$ is
the set of fractional ideals in $k$ that become principal over
$K$. It follows that $H^1(\gal,\oink_K^*)$ is in correspondence
with a subgroup of the ideal group of $k$.
\end{ex}

\begin{ex}
Let $R$ be a ring with a $\gal$-action. Then there is a natural action of $\gal$
on $A=\matrici_n(R)$ such that the matrices $E_{i,j}$ with a $1$ in the intersection of
the $i$-th row and the $j$-th column and $0$ everywhere else are invariant.

Any $A$-module $B$ has a decomposition of the form
$B=\oplus_{i=1}^n E_{i,i}B$. Furthermore, we claim that
$E_{i,j}(E_{k,k}B)=\delta_{j,k}E_{i,i}$. In fact, the case $j\neq k$ is trivial
and the case $j=k$ follows from the contentions
$E_{i,j}B=E_{i,i}(E_{i,j}B)\subseteq E_{i,i}B$ and
$E_{i,i}B=E_{i,j}(E_{j,i}B)\subseteq E_{i,j}B$.
Conversely, if $B'$ is an $R$-module of rank $n$, then $B=\bigoplus_{i=1}^\infty B_i$
where $B_i\cong B'$ has a natural $A$-module structure satisfying
$E_{i,j}B_k=\delta_{j,k}B_i$.
It follows that
The cohomology set $H^1(\gal,GL_n(R))$ classifies free left
$(R,\gal)$-modules of rank $n$.
\end{ex}

\section{Invariant generators of vector spaces}

In this section, $K/k$ is a finite Galois filed extension.
Also, $V_K$ denotes a finite dimensional vector space.
We assume that the Galois group $\gal$ acts on $V_K$
in such a way that $(\lambda v)^\sigma=\lambda^\sigma v^\sigma$ for any
$\sigma\in\gal$, $v\in V_K$, and $\la\in K$.

\begin{lem}
The maps $h\mapsto h^\sigma$, are linearly independent.
\end{lem}

\subparagraph{Proof.}
Assume $\sum_{\sigma\in\gal}\alpha_\sigma h^\sigma=0$ for all $h\in K$.
Since any finite separable extension is simple, we can assume $K=k(\omega)$.
In particular, $\sum_{\sigma\in\gal}\alpha_\sigma \omega^{k\sigma}=0$.
However, the matrix with entries $\omega^{k\sigma}$ is a Vandermonde matrix with
non-zero determinant, whence $\alpha_\sigma=0$ for all $\sigma$.
\qed

\begin{prop}\label{invariantbasis}
The space $V_K$ has a basis of $\gal$-invariant vectors.
\end{prop}

\subparagraph{Proof.}
Let $V_k$ be the subspace of invariant vectors.
It suffices to prove that $\dim_k(V_k)=\dim_K(V_K)$. Let $b:V_K\rightarrow V_k$
be the map $b(v)=\sum_{\sigma\in\gal}v^\sigma$. If $\dim_k(V_k)<\dim_K(V_K)$,
there exists a non-trivial lineal form $u$ such that $u(b(v))=0$ for all $v$.
For any $h\in K$ we have
\[
0=u(b(hv))=\sum_{\sigma\in\gal}h^\sigma u(v^\sigma).
\]
Since the functions $h\mapsto h\sigma$ are linearly independent, it follows that
$u(v^\sigma)=0$ for all $v$ and all $\sigma$. In particular $u=0$.
\qed

\begin{cor}
There exists a $k$-subspace $V_k$ of $V_K$, such that
$V_K\cong K\otimes_k V_k$.
\qed
\end{cor}

\begin{cor}[Hilbert's Theorem 90]\label{Hilbert90}
$H^1(\gal,\rm{GL}_n(K))=\uno$.\qed
\end{cor}

\section{Cohomology of the group of units of an algebra}

In this section, $K/k$ is a finite Galois extension
over an infinite field $k$.

\begin{prop}
\label{general}
For any finite dimensional algebra $A$, defined over $k$,
 and any algebraic extension
$K/k$, it holds that $H^1(\galkk,A^*_K)=\uno$.
\end{prop}

\subparagraph{Proof.}
Assume first that $K/k$ is finite.
Let $B_K$ be a free left $(A_K,\gal)$-module of rank $1$. Then $B_K$
is a finite dimensional vector space over $K$ satisfying the hypotheses
of proposition \ref{invariantbasis}. It follows that $B_K$ has a basis
 of invariant vectors $\{v_1,\dots,v_n\}$. In particular,
$B_K\cong K\otimes_k B_k$ for some $k$-vector space $B_k$.
The set of generators of $B_K$ is a Zariski open set and therefore
it must contain an element of $B_k$ (see some reference).
In the general case, if $\{L\}$ is the set of finite subextensions of $K/k$,
then
\[
H^1(\gal,A_K^*)=H^1(\lim_{\leftarrow}\gal_{L/k},\lim_{\rightarrow} A_L^*)=
\lim_{\rightarrow} H^1(\gal_{L/k},A_L^*)=\uno.\qed
\]

\begin{cor}
If $A_k$ is a $k$-algebra which is the direct limit of a family $\{B_k\}$ of finite
dimensional algebras then  $H^1(\galkk,A^*_K)=\uno$.
\end{cor}

\subparagraph{Proof.}
In fact,
\[
H^1(\gal,A_K^*)=H^1(\gal,\lim_{\rightarrow} B_K^*)=
\lim_{\rightarrow} H^1(\gal,B_K^*)=\uno.\qed
\]

\begin{cor}
If $A_k$ is a $k$-algebra such that every finite subset of $A_k$ generates a finite
dimensional subalgebra then  $H^1(\galkk,A^*_K)=\uno$.
\end{cor}

\subparagraph{Proof.}
In this case, $A_k$ is the direct limit of its finite dimensional subalgebras.
\qed

Since $\lineK \cong (\End_K(V))^*$, Hilbert's theorem 90 is a particular case
of proposition \ref{general}.
However, proposition \ref{general} has many other aplications,
as can be seen in these notes.

Let $K/k$ be any field, then we have an exact sequence:
\[
0 \vaa \specialK \vaa \lineK \vaa K^* \vaa 0
\]
which gives a long exact sequence in cohomology:
\[
\linek\rightarrow k^*\rightarrow H^1(\gal_{K/k},\specialK)\rightarrow
H^1(\gal_{K/k},\lineK)=1.
\]
as the determinant map is always surjective, this proves:

\begin{prop}
\label{special}
for any field extension $K/k$ we have:
\[
 H^1(\gal_{K/k},\specialK)=\{1\}.
\]
\end{prop}

More generally, if $A$ is a finite dimensional central simple algebra split by $K/k$
and $N_E:A_E \vaa E^*$ is the reduced norm, we have a sequence:
\[
A^*_k \stackrel{N_k}\vaa k^* \vaa H^1(\gal,\rm{S}A_K^*) \vaa H^1(\gal, A_K^*)=\{1\},
\]
where $\rm{S}A_K^*=ker(N_K)$ , therefore:
\[
H^1(\gal,\rm{S}A_K^*)=k^*/N_k(A^*_k).
\]

\begin{ex}
$A_E=\big(\frac{-1,-1}{E}\big)$, $A_\reali = \cuaternioni$ , $N(A_\reali)=\reali^+$,
and $A_\compleji=\matrici_2(\compleji)$, hence
$H^1(\gal,S\cuaternioni_\compleji^*)=\{[1],[-1]\}$.
Same result applies to any matrix algebra over $\cuaternioni$.

\end{ex}

\begin{prop}
The multiplicative and additive groups $\grupi_m$
and $\grupi_a$ have trivial $H^1$.
\end{prop}

\subparagraph{Proof.}
It is an immediate aplication of proposition \ref{general} that
$H^1(\gal,\grupi_m)=\uno$.

Let $A_k=k[x]/(x^2)$. $A$ is a local algebra with maximal ideal $I=(x)$.
Let $\mathcal{U}=\{1+y|y\in I \}$. Then, $\mathcal{U}\cong\grupi_a$.
Therefore, there exists an exact sequence
\[
\uno\vaa\grupi_a\vaa A^*\vaa\grupi_m\vaa\uno.
\]
It follows that $H^1(\gal,\grupi_a)\cong\coker(A_k^*\rightarrow k^*)=\uno$.
\qed

In fact, it holds that $H^i(\gal,\grupi_a)=\uno$ for all $i>0$. One way to
prove this is to see that $K$ is an induced $\gal$-module so that Shapiro's
lemma applies (\cite{kenneth}, p.73). If $k$ has characteristic $0$, an alternative
proof follows
from the fact that $H^i(\gal,A)$ is annihilated by $|\gal|$ for all $\gal$-module
$A$, while the map $\lambda\rightarrow n\lambda$ is an isomorphism for all $n$
(\cite{kenneth}, p. 84).

\begin{prop}
If $V_K$ is a finite dimensional vector space over $K$ provided with a Galois
action, then $H^1(\gal,V_K)=\uno$.
\end{prop}

\subparagraph{Proof.}
Let $\{v_1,\dots,v_n\}$ be a basis of invariant vectors of $V_K$, and let
$W=\spaan\{v_1,\dots,v_{n-1}\}$. There exists an exact sequence
\[
\uno\vaa W\vaa V\vaa\grupi_a\vaa\uno.
\]
Hence, the result follows by induction.

\begin{prop}
If $V_k$ is an arbitrary vector space over $k$, then $H^1(\gal,V_K)=\uno$.
\end{prop}

\subparagraph{Proof.}
This follows from the previous result since $V_k$ is a direct limit of its finite dimensional
subspaces $\{W_k\}$ and therefore
\[
H^1(\gal,V_K))=H^1(\gal,\lim_{\rightarrow} W_K)=\lim_{\rightarrow} H^1(\gal,W_K)=\uno.\qed
\]

\begin{ex}
Let $G$ be a $\gal$-group and let $K$ be a field with a faithfull $\gal$-action. Let $F=K^{\gal}$. It follows
from Galois theory that the extension $K/F$ is Galois and $\mathrm{Gal}(K/F)\cong\gal$.
Let $A_K=K[G]$ be the group algebra. Let $A_F=A_K^\gal$. It follows from proposition \ref{invariantbasis}
that $\dim_F(A_F)=\dim_K(A_K)$. In particular, $A_K$ is actually obtained from $A_F$ by extension of scalars and
proposition \ref{general} applies.
A basis $S$ of $A_K$ satisfying both $S=S^\sigma$ for $\sigma\in\gal$, and $GS=S$, is a principal homogeneous
$G$-space.
If a principal homogeneous $G$-space in isomorphic to some basis $S$ as above we say that it is represented in $A_K$.
Let $T$ be the set of bases $S$ of $A_K$ satisfying $GS=S$. then $A_K^*$ acts on $T$ by $S\mapsto Sa$ for $a\in A_K^*$.
Invariants elements of $T$ are principal homogeneous $G$-spaces represented in $A_K^*$.
By propositions \ref{classificationlaw} and \ref{general}, it follows that $T^\gal/A_F^*$ is in correspondence with
$H^1(\gal,G)$. It follows that every principal homogeneous $G$-space is represented in $A_K$ and if two bases $S$ and
$S'$ that are principal homogeneous $G$-spaces are isomorphic as such, then there exists $a\in A_F^*$ such that
$S'=Sa$.
\end{ex}

\chapter{Algebraic aplications of cohomology}

\section{Tensors and $K/k$-forms}\label{kkforms}\index{tensors}

By a tensor of type $(l,m)$ on $V$,
we mean an $\Omega$-linear map
$\tau:V^{\otimes l}\vaa V^{\otimes m}$,
 where
\[
V^{\otimes r}=\bigotimes_{i=1}^rV
\textnormal{ for } r\geq1,\ V^{\otimes 0}=\Omega.
\]
$\tau$ is said to be defined over $k$, if $\tau(V_k^{\otimes l})\subseteq
V_k^{\otimes m}$. \emph{All tensors mentioned in this work are assumed to
be defined over $k$}.
 $\lineomega$ acts on the set of tensors of type
$(l,m)$ by $g(\tau)=g^{\otimes m}\circ\tau\circ(g^{\otimes l})^{-1}$.
It makes sense, therefore, to speak about the stabilizer of a tensor.

Let $I$ be any set. By an $I$-family of tensors, we mean a map that
associates, to each element $i\in I$, a tensor $t_i$ of type $(n_i,m_i)$.
$\lineomega$ acts on the set of all $I$-families by acting in each
coordinate.
In all that follows, we say a family instead of an $I$-family unless
the set of indices needs to be made explicit.
Let $\tn$ be a family of tensors and $H=Stab_{\lineomega}(\tn)$.
Then, $H$ is a linear algebraic group.

If $K/k$ is a Galois extension with Galois group $\gal$,
 we get an exact
 sequence
\[
\uno\vaa H_K\vaa\lineK\vaa X_K\vaa\uno,
\]
where $X_K$ is the $\lineK$-orbit of $\tn$.
It follows from (\ref{clasificador}),
and example \ref{Hilbert90}, that
$X_K^\gal/\linek\cong H^1(\gal,H_K)$.
The elements of  $X_K^\gal/\linek$ can be thought of as
isomorphism classes of pairs $(V'_k,\tn')$ that become isomorphic to
$(V_k,\tn)$ when extended to $K$. These classes are usually called
$K/k$-forms of $(V,\tn)$, or just $k$-forms if $K=\bark$.
\index{kkforms@$K/k$-forms}
Observe that two vector spaces of the same dimension are isomorphic,
so we can always assume that the vector space $V$, in which all
tensors are defined, is fixed.

\begin{dfn}
We call a pair $(V,\tn)$, where $\tn$ is a family of tensors on $V$,
a \emph{space with tensors}, or simply a \emph{space}.
By abuse of language, we identify $(V,\tn)$ and $(V,\tn')$ whenever
$\tn,\tn'$ are in the same $\linek$-orbit,
i.e., if they correspond to the same $K/k$-form.
We say that $(V,\tn')$ is a $K/k$-form of $(V,\tn)$,
if $\tn$ and $\tn'$ are in the same $\lineK$ orbit.
\end{dfn}

\begin{ex}
Let $Q$ be a non-singular quadratic form on the space $V$. Then,
$\orthomega=Stab_{\lineomega}(Q)$.
Equivalence classes of non-singular
 quadratic forms on $V_k$ are classified by
$H^1(\gal,\orthkbar)$. A space, in this case, is what is usually called a
quadratic space.
\end{ex}

\section{Semi-simple abelian algebras}

An abelian algebra is semisimple if it has no nontrivial nilpotent elements.
In this section, let $L$ be an abelian semi-simple algebra
defined over a number field $k$. Then
 $L$ is a $k$-form of the trivial semi-simple algebra
\[
A^{(m)}=
\underbrace{k\oplus k\oplus \dots\oplus k}_{m\textnormal{ \scriptsize{times}}},
\]
whose group of automorfisms equals $S_m$, the simetric group on $m$ symbols.
It follows that the set of isomorfism clases of semi-simple abelian algebras
of dimension $m$ over $k$ is in one-to-one correspondence with the cohomogical
set $H^1(\gal,S_m)\cong Hom(\gal,S_m)/\sim$, where $\phi\sim\psi$ means that
there exists $\sigma\in S_m$ such that $\psi(g)=\sigma\phi(g)\sigma^{-1}$
for any $g\in\gal$.

Let $\psi:\gal\rightarrow S_m$ be one map in the conjugacy class corresponding to
$L$. then some properties of the algebra $L$ can be translated into properties of
the map $\psi$

Acording to the general theory, the algebra $L$ can be defined as the set of invariant
points of the corresponding twisted action, i.e.,
\[
L=\{l\in\bar{k}^m|\psi(\sigma)l^\sigma=l\ \forall\sigma\in\gal\}.
\]
Next we describe some properties of $L$ in terms of $\Psi=\im(\psi)$.

\begin{prop}
$L$ is a field iff and only if $\Psi$ acts transitively on the set
$\{1,\dots,m\}$.
\end{prop}

\subparagraph{Proof.}
The algebra $L$ is a field if it contains a non trivial projection $P$.
Let $P_t\in\bar{k}^m$ be the proyection in the $t$-th factor.
If $S\subseteq\{1,\dots,m\}$, define $P_S=\sum_{t\in S}P_t$.
Any proyection $P\in\bar{k}^m$ is of the form $P_S$ for some subset $S$.
Since all projections are fixed by the non-twisted action, $P_S\in L$ if and
only if $\psi(\sigma)P_S=P_S$ for all $\sigma\in\gal$. We have $\psi(\sigma)P_m=
P_{\psi(\sigma)(m)}$. Therefore, a non-trivial projection exists if and only if
$\Psi$ is not transitive.
\qed

More precisely, if $O\subseteq\{1,\dots,m\}$, then $P_O\in L$ if and only if
$O$ is invariant under $\Psi$, i.e., is a union of orbits.
Any element $l$ of $L$ has the form $l=(l_1,\dots,l_m)$ where the elements in
$l_1,\dots,l_m$ corresponding to elements in the same orbit form a complete set
of conjugates under the action of the Galois group, which acts on them by
$l_j^\sigma=l_{\psi^{-1}(j)}$. In particular, if  $L$ is a field,
$\{l_1,\dots,l_m\}$ is a complete set of conjugates.

\begin{prop}
Let $\hal$ be the subset of $\gal$ that fixes a subfield $L'$ of $\bar{k}$
isomorphic to $L$. Then
\[
\ker(\psi)=\bigcap_{\sigma\in\gal}\sigma\hal\sigma^{-1}.
\]
\end{prop}

The description of $L$ given earlier, implies that $L'$ can be assumed to
 be the image of $L$ under the projection on the first factor.
 the group $\ker(\psi)$ is the subgroup of $\gal$ that fixes $L$ pointwise
(in the non-twisted action).
Every conjugate of $\hal$ is the stabilizer of the image of $L$ under some projection.
\qed

\begin{cor}
If $L$ is a Galois extension of $k$, then $Gal(L/k)\cong\Psi$.\qed
\end{cor}

\begin{cor}
$L$ is a cyclic extension of $k$ if and only if $\im\psi$ is generated by an $m$-cycle.
\end{cor}

\subparagraph{Proof.}
If $L/k$ is a cyclic extension, then $\im(\psi)$ is cyclic and transitive.
\qed

\begin{ex}
The exact sequence
\[
\uno\vaa A_m\vaa S_m\vaa\mu_2\vaa\uno
\]
defines a cohomology map $d:H^1(\gal,S_m)\vaa H^1(\gal,\mu_2)\cong k^*/(k^*)^2$,
called the \emph{discriminant}, whose kernel is $H^1(\gal,A_m)$.
In other words, $L$ has discriminant $1$ if and only if $\Psi\subseteq A_m$.
It follows that if $L$ is Galois and $n$ is odd then the discriminant of $L$ is
$1$. The converse is true for $n=3$.
\end{ex}

\section{finite dimensional abelian algebras with nilpotent elements}

Let $L_k$ be an arbitrary finite dimensional abelian algebra over $k$, and let
$K$ be an algebraic closure of $k$.

\begin{prop}
There are no non-trivial $K/k$-forms of $k[x]/(x^n)$.
\end{prop}

\subparagraph{Proof.} Set $A_k=k[x]/(x^n)$. Then $A_K=K\oplus
I_K$, where $I_K$ is the ideal generated by $x$. Similarly
$A_k=k\oplus I_k$. Observe that $I_K^{n-1}\neq0$ and $I_K^n=0$. We
use induction on $n$. If $n=1$ there is nothing to prove. If
$n\geq2$, then $G_n=\Aut(A)$ acts on $A/I^{n-1}$. Let $G'$ be the
kernel of this action. Then $G_n/G'\cong G_{n-1}$, so that by
induction hipothesis $H^1(\gal,G/G')=\uno$.

Let $g\in G'$, and set $g(u)=u+e(u)$ with $e(u)\in I^m$.
Then $g(1)=1$, so that $e(1)=0$.
Also, for $uv\in I^2$,
$g(uv)=(u+e(u))(v+e(v))=uv$, so $e(uv)=0$.
It follows that $G'\cong\Hom(I/I^2,I^m)$. In particular, it is a vector space,
so it is acyclic.
\qed

A similar argument shows the following result:

\begin{prop}
There are no non-trivial $K/k$-forms of
$$k[x_1,\dots,x_n]/(x_1,\dots,x_n)^m.$$
\end{prop}

In general, for an algebra $A_k$ without projectors on the algebraic closure
$\bar{k}$ it is always true that
$A_K=K\oplus I_K$, ehere $I_K$ is the nilradical of $A_K$.
We can also define $G_n=\Aut(A_k/I_k^{n+1})$. However,
in general $G_n/G'$ is only a subgroup of $G_{n-1}$, and
it is not always true that $H^1(\gal,\Aut(A_K))=\uno$.

The following example will make this clear:

Let $A_k=k[x,y]/\Big((x,y)^3+(x^2-y^2)\Big)$. Then if $I_k$ is the
image of the ideal $(x,y)$, then $I^3=0$, but $I^2\neq 0$. Let $G$
be the automorphism group of $A$ and let $G'$ be the subgroup of
automorphisms of $A$ that induce the trivial automorphism of
$A/I^2$. Then $G/G'$ is contained in the group of automorphisms of
$A/I^2$, i.e., the group $\mathrm{GL}(V)$, where $V$ is the vector
space generated by $x$ and $y$. An element $g\in\mathrm{GL}(V)$ is
in $G/G'$ if and only if it fixes the ideal $(x,y)^3+(x^2-y^2)$.
It follows that there is a short exact sequence
\[
K^*\hookrightarrow G_K/G'_K\twoheadrightarrow\oink_K(q),
\]
where $q$ is the quadratic form $x^2-y^2$ and $\oink(q)$ its
orthogonal group. Since $K^*$ is acyclic, it follows that
$H^1(\gal,G_K)$ equals $H^1(\gal,\oink_K(q))$. In other words, the
$K/k$-forms of $A_k$ are the algebras of the form
$k[x,y]/\Big((x,y)^3+(q'(x,y))\Big)$, where $q'$ is a quadratic
form.

\section{Skolem-Noether theorem}

Let $\alge$ be a central simple algebra defined over $k$.
Let $L$ be a maximal semisimple Abelian subalgebra defined over $k$.
$\alge^*$ acts on the set of maximal semisimple
Abelian subalgebras by conjugation.
It is a trivial exercise in linear algebra to prove the
transitivity of this action over an algebraically closed field.
Let $G$ be the stabilizer of $L$.
It follows from (\ref{clasificador}), and proposition \ref{general},
that the set of conjugacy classes of
maximal Abelian subalgebras that are defined over $k$
is parametrized by $H^1(\gal,G)$.
Observe that, over the algebraic closure,
any automorphism of $L$ arises from a conjugation.
 The short exact sequence
\[
\uno\vaa L^*\vaa G \vaa\Aut
(L)\vaa\uno,
\]
where $\Aut(L)$ is the set of automorphisms of $L$ as an
$\Omega$-algebra, gives
\[
1= H^1(\gal,L^*)\vaa H^1(\gal,G)\vaa H^1(\gal,\Aut(L)).
\]
It follows,
since $H^1(\gal,\Aut(L))$ classifies isomorphism classes of
semisimple Abelian algebras,
that any two isomorphic algebras are conjugate.

\section{Finite subgroups in proyective groups of algebras}

Let $k$ be a field and let $\alge_k$ be a finite dimensional
$k$-algebra. Let $\Gamma_0$ be a finite subgroup of
$\alge_k^*/k^*$. The group $\Gamma_0$ can be regarded as a
subgroup of $\alge_K^*/K^*$ for any field extension $K/k$. In this
section we describe a cohomology set that classifies finite
subgroups $\Gamma$ of $\alge_k^*/k^*$ that become conjugate to
$\Gamma_0$ over some separable algebraic extension $K/k$. Let
$C_K$ be the centralizer of $\Gamma_0$ in $\alge_K^*/K^*$ and let
$W$ be the group of automorphisms of $\Gamma_0$ that can be
realized as conjugations from elements in $\alge_K^*/K^*$. In this
section we prove the following result:

\begin{prop}\label{thm1}
Let $K/k$ be a Galois extension. There exists a natural action of
$W$ on the cohomology set $H^1(C_K)=H^1(K/k,C_K)$. The set of
conjugacy classes of finite subgroups of $\alge_k^*/k^*$ that
become conjugate over $K$ to $\Gamma_0$ is in one-to-one
correspondence with the set of orbits $H^1(C_K)/W$.
\end{prop}

\subparagraph{Proof.} The group $\alge_K^*$ act on the set of
finite subgroups by conjugation. The stabilizer of $\Gamma_0$ is
the preimage $N_K\subseteq\alge_K^*$ of the normalizer of
$\Gamma_0$ in $\alge_K/K^*$. It follows from Proposition
\ref{classificationlaw} that, if $X$ is the set of finite
subgroups $\Gamma$ of $\alge_K^*/K^*$ that are conjugate to
$\Gamma_0$, i.e., the $\alge^*_K$-orbit of $\Gamma_0$, then
$X^\gal/\alge^*_k\cong\ker[H^1(N_K)\rightarrow H^1(\alge_K^*)]$.
Since $H^1(\alge_K^*)=\uno$ (\cite{kneser}, p. 16), it follows
that $H^1(N_K)$ is in correspondence with the set of conjugacy
classes under $\alge^*_k$ of $\gal$-invariant finite groups
$\Gamma$ that are $\alge^*_K$-conjugate to $\Gamma_0$. The initial
group $\Gamma_0$ corresponds to the distinguished element of
$H^1(N_K)$. There exists a short exact sequence
$C_K\hookrightarrow N_K \twoheadrightarrow W$, where $W$ is a
subgroup of $\Aut(\Gamma_0)$. Notice that $\gal$ acts trivially on
$\Aut(\Gamma_0)$, and hence also  on $W$. It follows from
Propositions 38 and 39 in (\cite{Serre}, p. 49), that
$H^1(C_K)/W\cong\ker[H^1(N_K)\stackrel{\pi}{\rightarrow}H^1(W)]$.
Now, since the action of $\gal$ on $W$ is trivial, the set
$H^1(W)$ is identified with the set of conjugacy classes of
homomorphisms from $\gal$ to $W$. Under this identification, the
map $\pi$ sends a cocycle $\alpha\in H^1(N_K)$ to a map
$\phi_\alpha$, where $\phi_\alpha(\sigma)\in W$ acts on $\Gamma_0$
as conjugation by $\alpha_\sigma$. If $\alpha$ is the cocycle
corresponding to a group $\Gamma=a\Gamma_0a^{-1}$, then
$\alpha_\sigma=a^{-1}a^\sigma$. Now,
\[
(a\gamma a^{-1})^\sigma
=a\alpha_\sigma\gamma\alpha_\sigma^{-1}a^{-1}\qquad\forall\gamma\in\Gamma_0
\subseteq\alge_k.
\]
It follows that the kernel of the map $\phi_\alpha$ is the
subgroup of $\gal$ corresponding to
 the Galois extension $k(\Gamma)/k$ generated by
the coordinates of the elements of $\Gamma$. Hence $\phi_\alpha$
is trivial if and only if $k(\Gamma)=k$. The result follows.\qed

 We give two applications of this result:

\begin{cor}\label{thm2}
Let $k$ be a field whose characteristic does not divide $n$. If
$\alge_k$ is a central division algebra over $k$ of dimension
$n^2$, and if $\Gamma_0$ contains a subgroup $\Omega\cong
\enteri/n\enteri\times\enteri/n\enteri$ that intersects trivially
the center $Z(\Gamma_0)$ of $\Gamma_0$, then every finite group
$\Gamma$ of $\alge_k^*/k^*$ that is conjugate to $\Gamma_0$ over
$K$ is conjugate to $\Gamma_0$ over $k$.
\end{cor}

\subparagraph{Proof.} Without loss of generality we assume $K/k$
is a Galois extension. we claim that if $\Omega$ is a subgroup of
$\alge_k^*/k^*$ isomorphic to $\enteri/n\enteri\times
\enteri/n\enteri$, then $\Omega$ is its own centralizer in
$\alge_K^*/K^*$. The group $\Omega$ is generated by $\rho(x)$ and
$\rho(y)$, where $\rho:\alge_K^*\rightarrow:\alge_K^*/K^*$ is the
canonical projection, and the elements $x$ and $y$ satisfy the
identities $x^n=a$, $y^n=b$, and $xy=\eta yx$ for some
$a,b,\eta\in k$. By Kummer's Theory, $y^t\in k(x)$ implies that
$n$ divides $t$. In particular, $\eta$ is a primitive $n$-th root
of unity. By considering the eigenvectors of $u\mapsto xux^{-1}$,
it follows that $x$ and $y$ generates $\alge_K$ as a $K$-algebra.
Any element centralizing $\rho(x)$ and $\rho(y)$ must be of the
form $\rho(z)$ where $z^{-1}xz=\tau x$ and $z^{-1}yz=\kappa y$.
Comparing eigenvalues of the functions $u\mapsto xu$ and $u\mapsto
(z^{-1}xz)u$ we see that $\tau=\eta^t$, and similarly
$\kappa=\eta^s$. It follows that $zy^{-t}x^{-s}$ centralizes both
$x$ and $y$, so that $\rho(z)=\rho(y)^t\rho(x)^s$. The condition
that $\Omega$ meets trivially the center of $\Gamma_0$ implies
that $C_K=K^*$, hence $H^1(C_K)=\{1\}$ (\cite{kneser}, p. 16). The
result follows by Theorem \ref{thm1}. \qed

This result applies to groups $\Gamma_0$ containing a copy of
$A_{2n}$ for $n^2=\dim_k\alge_k>1$. If $n=2$, it applies to groups
isomorphic to $A_4$, $S_4$ or $A_5$.

\begin{cor}\label{thm3}
Let $k$ be a field of characteristic not equal to $p$, where $p$
is a prime. If $\alge_k$ is a central division algebra over $k$ of
dimension $p^2$, and if $\Gamma_0$ is isomorphic to
$\enteri/m\enteri$ with $m\neq p$, then every finite group
$\Gamma$ of $\alge_k^*/k^*$ that is conjugate to $\Gamma_0$ over
$K$ is conjugate to $\Gamma_0$ over $k$.
\end{cor}

\subparagraph{Proof.} Without loss of generality we assume $K/k$
is a Galois extension. The group $\Gamma_0$ is generated by an
element $\rho(x)$, where $x\in\alge_k$ satisfies $x^m\in k$ and
$\rho$ is as in the proof of Corollary \ref{thm2}. Since $\alge_k$
is a division algebra, for any $d<m$ the subalgebra $k(x^d)$ is
maximal abelian in $\alge_k$, i.e., $[k(x^d):k]=p$. In particular,
$k(x^d)=k(x)$. Furthermore, we have $K(x)^*\subseteq C_K$.
 Then any element $z\in\alge_K$
such that $\rho(z)$ centralizes $\rho(x)$ must satisfy
$z^{-1}xz=\tau x$ where $\tau$ is an $m$-th root of unity. Since
$\tau x$ and $x$ have the same eigenvalues, multiplication by
$\tau$ must permute the eigenvalues of $x$. It follows that the
order of $\tau$ must divide $p!$. Assume first that $m$ is not a
power of $p$. Replacing $x$ by some power if needed we might
assume that $m$ is a prime. Since $[k(x):k]=p$, we must have
$m>p$. We conclude that $z\in K(x)^*$. Assume next $m=p^r$ with
$r>1$. The same argument as above shows that $\tau$ must be a
$p$-th root of unity. Hence $z$ commutes with $x^p$. Since
$k(x^p)=k(x)$, we still have $z\in K(x)^*$. As $H^1(K(x)^*)=1$
(\cite{kneser}, p. 16), the result follows. \qed

\chapter{Lattices and cohomology}\label{chapter5}

\section{Basic results}

Let $k$ be a local or number field, $K/k$ a Galois extension,
 $G\subseteq\lineomega$
 an algebraic group
defined over $k$, $\la_k$ a lattice on $V_k$, $L_K$  a $\gal$-invariant
lattice on $V_K$. Let $\gal=\galkk$.

\begin{prop}\label{result1}
If there is an element $\varphi\in G_K$ such that
$\varphi(L_K)=\Lambda_K$, then
$a_\sigma=\varphi^\sigma\varphi^{-1}$ is a well defined element of
$H^1(\gal,G_K^\la)$. It is independent of the choice of
$\varphi$ and depends only on the orbit of $L_K$ under $G_k$.
The correspondence assigning, to every such
$G_k$-orbit of $\oink_K$-lattices,
 an equivalence class of cocycles,  is an injection.
The image of this map is
\[
\ker(H^1(\gal,G_K^{\la})\stackrel{i_*}{\longrightarrow}H^1(\gal,G_K)),
\]
where $i$ is the inclusion.
\end{prop}

\paragraph{Proof.}
$G_K$ acts on the set of $\oink_K$-lattices in $V_K$.
Let $X$ be the orbit of $\Lambda_K$.
 We have an exact sequence
\[
\uno\vaa G_K^{\la}\vaa G_K \vaa X\vaa\uno.
\]
Hence, by
(\ref{clasificador}),
we get
$X^\gal/G_k \cong \ \ker(H^1(\gal,G_K^{\la}) \vaa
H^1(\gal,G_K))$.
$\Box$

\begin{ex} Using the fact that $H^1(\gal,\lineK)=\{1\}$,
we obtain that the set of $\linek$-orbits of $\gal$-invariant $\oink_K$-lattices
isomorphic to $\la_K$ is in correspondence with $H^1(\gal,\lineKlam)$.
\end{ex}

If $G$ is defined as the stabilizer of a family of
 tensors, e.g., the unitary group of a hermitian form or the
automorphism group of an algebra, we get
a more precise result.

Recall that in section \ref{kkforms} we identified
$K/k$-forms of $(V,\tn)$ with the corresponding
$\linek$-orbits of families of tensors.

\begin{dfn}
Let $(V,\tn)$ be a space.
 A lattice in $(V_K,\tn)$
is a pair $(\la_K,\tn)$, where $\la_K$ is a lattice in $V_K$.
$\lineK$ acts on the set of pairs $(\la_K,\tn')$,
for all families of tensors $\tn'$, by acting on each component.
 Two lattices $(\la_K,\tn)$, $(L_K,\tn')$
are said to be \emph{in the same space} if $\tn,\tn'$
 are in the same $\linek$-orbit.
\end{dfn}

\begin{prop}\label{section3.4}
 Assume that $G$ is the
stabilizer of a family of tensors $\tn$ on $V$.
The set $H^1(\gal,G_K^{\la})$ is
in one-to-one correspondence
with the set of $G_k$-orbits of $\gal$-invariant
$\oink_K$-lattices in the same $G_K$-orbit,
in all spaces that are $K/k$-forms of $(V,\tn)$.
The kernel of the map
\[
H^1(\gal,G_K^\la)\stackrel{i_*}{\longrightarrow}H^1(\gal,G_K),
\]
where $i$ is the inclusion, corresponds to
the subset of orbits of lattices that
are in the
same space as $\la_K$.
\end{prop}

\paragraph{Proof.}
We have an action of $\lineK$ on the set of all pairs $(L_K,\tn')$,
where $L_K$ is a lattice and $\tn'$ a family of tensors with a
fixed index set. If $T$ is the orbit of
$(\la_K,\tn)$, we have a sequence
\[
\uno\vaa G_K^\la\vaa\lineK\vaa T\vaa\uno,
\]
and the same argument as before applies. Last statement
 follows from the fact that spaces $(V_K,\tn')$ are classified by
$H^1(\galkk,G_K)$, (see section \ref{kkforms} or \cite{kneser}, p. 15).
$\Box$

\begin{rmk}
Recall that $\la_K=\la_k\otimes_{\oink_k}\oink_K$.
If $L_K$ is in the same $G_k$-orbit as $\la_K$, $L_K=L_k
\otimes_{\oink_k}\oink_K$, since $G_k$ also acts on $V_k$.
Recall that we defined the cocycle corresponding to $L$
by the formula $a_\sigma=\phi^\sigma\phi^{-1}$ (see prop.
\ref{result1}). This definition does not depend on $G$,
as long as $\phi\in G$.
It follows that the set of $G_k$-orbits of
lattices in $V_k$ that are isomorphic as $\oink_k$-modules,
and whose extensions to $K$ are in the same $G_K$ orbit,
corresponds to
\begin{equation}\label{kernelintersection}
\ker\left(H^1(\gal,G_K^\la)\vaa H^1(\gal,G_K)\times
H^1(\gal,\lineKlam)\right).
\end{equation}
In the case that $G$ is the stabilizer of a family of tensors,
\[
\ker(H^1(\gal,G_K^\la)\vaa H^1(\gal,\lineKlam))
\]
corresponds to the set of $G_k$-orbits of
such lattices in all spaces that are
$K/k$-forms of $(V,\tn)$.
\end{rmk}

\begin{ex}
If $\la_k$ is free, (\ref{kernelintersection})
 corresponds to the set of $G_k$-orbits of free lattices
on $V_k$, whose extensions to $K$
are in the same $G_K$-orbit.
\end{ex}

\begin{dfn}
We say that an $\oink_K$-lattice $\la_K$ is
\emph{defined over $k$}, if
$\la_K\cong\oink_K\otimes_{\oink_k}\la_k$ for some $\la_k$.
We say that $\la_K$ is a \emph{$k$-free lattice}, if $\la_k$ is free.
\end{dfn}

Assume first that $G$ is the stabilizer of a family of tensors.

\begin{dfn}
Let $a\in H^1(\gal,\gKla)$. We say that $a$ is defined over
$k$, $k$-free or
in $(V,\tn)$ if some (hence any),
lattice in the class corresponding to $a$ has this property. Define
\index{lk@$\lk,\lfree$, etc}
\[
\begin{array}{ccl}
\Lk G{K/k}{\la}&=&\{ a\in H^1(\gal,\gKla) | a \textnormal{
is defined over } k\},\\
\Lfree G{K/k}{\la}&=&\{a\in \Lk G{K/k}{\la}| a
\textnormal{  is $k$-free}\},\\
\Lvk G{K/k}{\la}&=&\{a \in  H^1(\gal,\gKla)| a
\textnormal{  is in } (V_K,\tn)\},\\
\Lkvk G{K/k}{\la}&=&\Lvk G{K/k}{\la}\cap \Lk G{K/k}{\la},\\
\Lfreevk G{K/k}{\la}&=&\Lvk G{K/k}{\la}\cap\Lfree G{K/k}{\la}.
\end{array}
\]
\end{dfn}

Let
\begin{equation}
\label{eq1}
F_1:  H^1(\gal,\gKla)\longrightarrow H^1(\gal,\gK),
\end{equation}
\begin{equation}
\label{eq2}
F_2:  H^1(\gal,\gKla)\longrightarrow H^1(\gal,\lineKlam),
\end{equation}
be the maps defined by the inclusions. Then, we have the following
proposition:

\begin{prop}\label{proposition7}
Assume that $\la_k$ is free.
The following identities hold:
\[
\begin{array}{ccl}
\Lvk G{K/k}{\la} &=& \ker F_1,\\
\Lfree G{K/k}{\la} &=& \ker F_2,\\
\Lfreevk G{K/k}{\la} &=& \ker F_1\cap\ker F_2.\Box
\end{array}
\]
\end{prop}

Later, we give a similar interpretation to $\lk$.

\begin{ex}
\[
\Lfree {\orthomega}{\bark/k}{\la}=
\ker(H^1(\gal,\orthi n{\bark}Q\la)\vaa H^1(\gal,\linekbarlam))
\]
is in correspondence with the set of isometry classes of
free quadratic lattices that become isometric to $\la_k$ over
some extension.
\end{ex}

\begin{rmk}
Notice that $\lvk,\lkvk,\lfreevk$ can be defined, even if $G$
is not the stabilizer of a family of tensors, as follows:
\[
\begin{array}{rcl}
\Lvk G{K/k}{\la}&=&\ker(H^1(\gal,\gKla)\vaa H^1(\gal,\gK)),\\
\Lkvk G{K/k}{\la}&=&\{a\in\Lvk G{K/k}{\la}| a \textnormal{ is
 defined over }k\} ,\\
\Lfreevk G{K/k}{\la}&=&\{a\in\Lvk G{K/k}{\la}| a \textnormal{ is
 free}\}.
\end{array}
\]
In this case, the first and last identities of proposition
\ref{proposition7} still hold. Notice that we can still interpret
 $\lvk$ as a set of equivalence classes of lattices,
because of proposition \ref{result1}.
\end{rmk}

\section{$H^1(\gal,U_K)$ and the ideal group}\label{units}

Let $k$ be a local or number field,
$K/k$ a finite Galois extension.
Let  $\gal=\galkk$. $U_K=\oink_K^*$ denotes the group of units of $\oink_k$.
\index{uk@$U_k$}

For any local or number field $E$, let $I_E$ be its group of fractional
ideals, $P_E$ the subgroup of principal fractional ideals.
There is a natural map
$\alpha:I_k\rightarrow I_K$ defined by
$\alpha(\ideala)=\ideala\tensor\oink_K$.  Clearly
$\alpha(P_k)\subseteq P_K$,
so we get a map $\alpha': I_k/P_k\rightarrow I_K/P_K$.

Apply the general theory to $\la_k=\oink_k$, $G=\grupi_m$,  $G_K^{\la}=U_K$.
Any $\lambda\in K^*$ acts by $\ideala\mapsto\lambda\ideala$, for $\ideala\in I_K$.
It follows that,
\[
H^1(\gal,U_K)\cong(P_K)^\gal/\alpha(P_k).
\]

\begin{rmk}
When this is done over the function field extension $L(X)/F(X)$,
one obtains $U_K=L^*$. Since $H^1(\gal,L^*)=\{1\}$, this implies
that every invariant ideal in $L[X]$ has a generator in $F[X]$. In
particular, if a polinomial $P$ of $F[X]$ is an $n$-power in
$L[X]$, by taking the principal ideal generated by the $n$-th
root, we obtain that $\lambda P$ is an $n$-th power in $F[X]$ for
some constant $\lambda$. Since powers of monic polinomials are
monic we conclude that a monic in $F[X]$ is an $n$-th power in
$L[X]$ if and only if it is an $n$-th power in $F[X]$.
\end{rmk}

Non-zero prime ideals of $\oink_K$ form a set of free generators for $I_K$
(see \cite{Lang}, p. 18).
Let $\ideala\in I_K$. We can write
\[
\ideala
=\prod_{\idealp\in \Pi(k)}(\prod_{\idealP | \idealp}
\idealP^{\beta(\idealP)}).
\]
 If $\ideala$ is $\gal$-invariant, all
the powers $\beta(\idealP)$ corresponding to prime
 divisors of the same prime of $k$ must be equal.  In other words:
\begin{equation}
\label{productus}
\ideala=\prod_{\idealp\in\Pi(k)}(\prod_{\idealP |
\idealp}\idealP)^{\beta(\idealp)},
\end{equation}
where $\beta(\idealp)$ is the common value of $\beta(\idealP)$ for all
$\idealP$ dividing $\idealp$.
This ideal is in $\alpha(I_k)$ if and only if the ramification degree
$e_\idealp$ divides
$\beta(\idealp)$ for all $\idealp$. Hence, we have an exact sequence
\[
0\longrightarrow \ker\alpha'\longrightarrow(P_K)^\gal/\alpha(P_k)
\longrightarrow\prod_{\idealp\in \Pi(k)}(\enteri/e_\idealp) ,
\]
where the image of the last map corresponds to those ideals of the form
 (\ref{productus}) that are principal in $K$.
The image of $\ker\alpha'$ is what we call
$\Lk G{K/k}{\la}$.
  In particular, since all ideals become principal in some
 extension, we can take a direct limit, to
obtain the long exact sequence:
\[
0\vaa I_k/P_k\vaa H^1(\galbarkk,U_{\bar{k}})
\vaa\prod_{\idealp\in \Pi(k)}(\Q/\enteri)\vaa 0.
\]
A refinement of this argument gives
\[
H^1(\galbarkk,U_{\bar{k}})\cong(I_k\bigotimes_\enteri\Q)/
(P_k\bigotimes_\enteri\enteri),\ \ \Lk G{K/k}{\la}=I_k/P_k.
\]

\section{Localization}

Recall remarks \ref{fixw} and \ref{galoisgroups}.

Assume $k$ is a number field.
There exist natural localization maps
\[
F_v:H^1(\gal,\gKla)\rightarrow H^1(\galv,\gKlav),
\]
defined by inclusion and restriction.
We define $\gKlav=G_{K_w}$ if $w$ is Archimedean.

\begin{lem}
\label{localkernelgglam}
Let $F_1:H^1(\gal,\gKla)\rightarrow H^1(\gal,G_K)$ be the map induced by
the inclusion. If the natural map
\[
\tau:H^1(\gal,G_K)\rightarrow \prod_{v\in\Pi(k)} H^1(\galv,G_{K_w})
\]
is injective, then $\ker F_1\supseteq \bigcap_v \ker F_v$.
\end{lem}

\subparagraph{Proof of lemma.}
Immediate from the following
commutative diagram:
\[
\xymatrix{
H^1(\gal,\gKla) \ar[r]^{F_1}\ar[d]^{\prod_vF_v}&
H^1(\gal,G_K)\ar[d]^{\tau} \\
\prod_v H^1(\galv,\gKlav)\ar[r] &
\prod_v H^1(\galv,G_{K_w}).\Box}
\]

\begin{rmk}
If the hypothesis of this lemma is satisfied, one says that $G$
 satisfies the Hasse principle over $k$.
\end{rmk}

\subparagraph{Characterisation of $\lk$.}

$\Lk G{K/k}{\la}$ is the set of equivalence classes
 of lattices defined over $k$
that become isomorphic over $K$.
A lattice $L_K$ is defined over $k$ if and only if it is generated by its
$k$-points, i.e.,
\[
L_K=\oink_K(L_K\cap V_k).
\]
This is a local property. On the other hand, for any local place $v$, all lattices defined over $k_v$ are $k_v$-free,
i.e.,
$$\Lk {\lineomega}{K_w/k_v}{\la}=\Lfree {\lineomega}{K_w/k_v}{\la}.$$
The following result is immediate from this observation.
\begin{prop}\label{lkernel}
\[
\Lk G{K/k}{\la}=\ker\left(H^1(\gal,\gKla)\vaa\prod_v
H^1(\galv,\lineKvlam)\right).\Box
\]
\end{prop}

\section{Genus and cohomology}\label{chapter6}

Assume that
In all of section \ref{chapter6}, $k$ is a number field.

\begin{dfn}
Let $F_v$ be the localization map. Define\index{cgenus@$C$-genus}
\[
\ckgen G{K/k}{\la}=\ker(\prod_vF_v).
\]
We call this set the \emph{cohomological genus} of $\la$ with respect to
$G$. \index{cohomological genus}
\end{dfn}

\begin{prop}\label{wcisc}
For any linear algebraic group $G$, it holds that
\[
\ckgen G{K/k}{\la}\subseteq
\Lk G{K/k}{\la}.
\]
\end{prop}

\paragraph{Proof.}
This follows from proposition \ref{lkernel} and the
commutative diagram
\[
\xymatrix{H^1(\gal,\gKla)\ar[d]\ar[rd]&\\
\prod_{v\in\Pi(k)}H^1(\galv,\gKlav)\ar[r]&
\prod_{v\in\Pi(k)}H^1(\galv,\lineKvlam).\Box}
\]

\begin{rmk}
Assume $G$ is the stabilizer of a family of tensors.
This result tells us that the cohomological genus corresponds to a
set of equivalence classes ol lattices defined over $k$.
In fact, $a\in\ckgen G{K/k}{\la}$ if and only if $a$ corresponds to a lattice,
in some $K/k$-form of $(V,\tn)$, that is in the same $G_{k_v}$-orbit,
at every place $v$.
\end{rmk}

\begin{dfn}\index{vcgenus@$VC$-genus}
We define the $VC$-genus of $\la_k$
 by the formula
\[
V\ckgen G{K/k}{\la}=
\ckgen G{K/k}{\la}\cap \Lvk G{K/k}{\la}.
\]
In other words, it is the kernel of the map
\begin{equation}\label{cohomologymap}
H^1(\gal,\gKla)\vaa H^1(\gal.\gK)\times\prod_{v\in\Pi(k)}
H^1(\galv,\gKlav).
\end{equation}
\end{dfn}

Let $G$ be an arbitrary linear algebraic group.
The $VC$-genus corresponds to a set of $G_k$-orbits
of lattices in $V_k$.
In fact, it corresponds to a subset of the set of double cosets
$\gk\backslash \gak/\gakla$, i.e., the \emph{genus} of $G$
(see \cite{Pla}, p. 440).
In particular, the following proposition holds.

\begin{prop}\label{prevstrongaproxy}
If $G$ has class number $1$ with respect to a lattice $\la_k$, then
(\ref{cohomologymap}) has trivial kernel for every Galois extension $K/k$
(compare with corollary 4 on p. 491 of \cite{Pla}).$\Box$
\end{prop}

This, in particular, applies to a group having absolute strong
approximation (see \cite{Pla}).
 However, we have a stronger result.

\begin{prop}
\label{strongaproxy}
If $G$ has absolute strong approximation over $k$, the map
(\ref{cohomologymap}) is injective.
\end{prop}

\subparagraph{Proof.}
Recall remark \ref{fixw}.

Let $M_K$,$L_K$ be two $\gal$-invariant $\oink_K$-lattices in
$V_K$, that are locally in the same $\gkv$-orbit for all $v$.
Then, we can choose
elements $\sigma_v\in\gkv$, such that $\sigma_v M_{K_w}=L_{K_w}$
for every place $v$, and $\sigma_v=1$ at all but finite places. Now, any global element
$\sigma$, close enough to $\sigma_v$ at all finite
places where $\sigma_v\neq1$,
and stabilizing $M_{K_w}=L_{K_w}$ at the remaining finite places,
satisfies $\sigma M_K=L_K$, as claimed.
$\Box$

The following result is just a restatement of lemma
\ref{localkernelgglam}.

\begin{prop}\label{vcisc}
If $G$ satisfies the Hasse principle over $k$, then
\[
V\ckgen G{K/k}{\la}=\ckgen G{K/k}{\la}.\Box
\]
\end{prop}

This result tells us that, in the presence of Hasse principle,
the cohomological genus corresponds to a subset of the genus
(compare with \cite{Rohlfs}, thm 3.3, p. 198).

\section{Spinor norm and genera}

 Let $G\subseteq GL(V)$ be a semi-simple group with universal cover
$\gtilde$ and fundamental group $\mu_n$. Let $K=\bark$.

The short exact sequence
\[
\uno\vaa \mu_n \vaa\gtilde_K\vaa G_K \vaa\uno,
\]
defines a map
$\theta:\gk\vaa H^1(\gal,F)=\ksobrekn$.

Let $\la_k$ be any lattice in $V_k$.
 The following proposition holds.

\begin{prop}
With the above notations,
$V\ckgen G{K/k}{\la}$ is in one-to-one correspondence with
the genus of $G$
\textnormal{(compare with theorem 8.13 in \cite{Pla}, p. 490)}.
\end{prop}

\subparagraph{Proof.}
It suffices to show that any two $G_k$-orbits in the same genus are
identified over some extension.
Without loss of generality, we assume $k$ is non-real. It suffices
to check that they are in the same spinor genus (see \cite{spinor}).
Spinor genera are classified by
\[
J_k/J_k^nk^*\Theta_{\ad}
(G^{\la}_{\ad_k}),
\]
where
 $\Theta_{\ad}(G^{\la}_{\ad_k})$
is the kernel of the local spinor norm
(see\footnote{ The case of an orthogonal group is already considered in \cite{Hsia99}.}
\cite{mithesis} or \cite{spinor}).
This is a finite set, and the representing adeles can be chosen to have
trivial coordinates at almost all places.
Therefore, it suffices to take an extension that contains the
$n$-roots of unity, and $n$ roots of a finite set of local elements.$\Box$

This result allows us to use cohomology to study the genus of any
Semisimple group.

\section{Determinant class of a lattice}\index{determinant class}

Let $[\ideala]$ be the $k^*$-orbit of the $\oink_K$-ideal $\ideala$.
Assume that
\[
\la_k=\underbrace{\oink_k\oplus\dots\oplus\oink_k}_{n\textnormal{
\scriptsize times}}.
\]

The map
$\dete_*:H^1(\gal,\lineKlam)\vaa H^1(\gal,U_K)$ is the map induced
 in cohomology by the determinant.
It is surjective, since $\dete$
has a right inverse.
However, in general it is not injective, as the example below shows.

\begin{dfn}
 Let $L_K$ be a $\gal$-invariant lattice in $V_K$, and let $a$ be the
cocycle class corresponding to the $\linek$-orbit of $L_K$.  We
define the determinant class of $L_K$, which we denote $\dete_*(L_K)$, by:
\[
\dete_*(L_K)=\dete_*(a)\in H^1(\gal,U_K)\cong P_K^\gal/\alpha(P_k),
\]
and we identify it with the corresponding ideal class.
\end{dfn}

\begin{ex}
Using the standard embedding $GL(V)\times GL(W)\vaa GL(V\oplus W)$,
it is easy to prove that  $\dete_*(\la_K\oplus L_K)=
\dete_*(\la_K)\dete_*(L_K)$. In particular, we obtain
that $\dete_*(\ideala_1\oplus\dots\oplus\ideala_n)=
[\ideala_1\dots\ideala_n]$.

Assume $k\subseteq K$ are local fields with maximal ideals $\idealp$,
$\idealP$. Assume that $\idealp\oink_K=\idealP^e$. Then,
\[
\dete_*(\underbrace{\idealP\oplus\dots\idealP}_e)=
[\idealP^e]=1=\dete_*(\underbrace{\oink_K\oplus\dots\oplus\oink_K}_e),
\]
but the latter lattice is defined over $k$
 and the first one is not.
\end{ex}

Let $\lk=\Lk {GL(V)}{K/k}{\la}$.
We have the following result:

\begin{lem}
\label{tiroliro}
$\lk\cap\ker(\dete_*)=\uno$.
\end{lem}

\paragraph{Proof of lemma.}

This follows from the fact
that all $k$-defined lattices are of the form
$\ideala_k\oplus\oink_k\oplus\dots\oplus\oink_k$
(see \cite{Om}, (81:5)).
It can also be proved by a diagram chasing argument.
$\Box$

Now observe that, for any algebraic group $G\subseteq GL(V)$,
we have $\Lk G{K/k}{\la}=i^{-1}_*(\lk)$,
where $i_*$ is the cohomology map induced by the inclusion.

\begin{prop}
\label{lkernel2}
If $G\subseteq\specialomega$, then $i_*^{-1}(\lk)=\ker(i_*)$.
\end{prop}

\paragraph{Proof of proposition.}
It is immediate from the commutative diagram
$$\xymatrix{H^1(\gal,\gKla)\ar[rr]\ar[dr]_{i_*}&
{}\ar@{}[d]|{\circlearrowright}&H^1(\gal,\specialKlam)\ar[dl]
\\{}&H^1(\gal,\lineKlam)\ar[d]^{\dete_*}&{}\\&H^1(\gal,U_K)&{}}$$
 that $\im(i_*)\subseteq\ker(\dete_*)$.
Now recall  lemma \ref{tiroliro}.$\Box$

In particular, such a group cannot identify a free lattice to a
non-free $k$-defined lattice over any extension, although it can
identify a free lattice to a non-$k$-defined lattice.

In this case, a description of $\lfree$ is equivalent to a description
of $\lk$, hence $\Lk {G}{K/k}{\la}$ can be described without resorting to
localization.

\section{Cohomology and representation}

Let $M$ be a sublattice of $\la$. Let $G_K^{\la,M}$ be the stabilizer
 of $M$ in $G_K^{\la}$. There exists a short exact sequence:
\[
\uno\vaa G_K^{\la,M}\vaa G_K^{\la}\vaa X\vaa\uno,
\]
where $X$ is the orbit of $M$ in the set of sublattices. Then
\[
X^\gal/G_k^{\la}\cong\ker(H^1(\gal,G_K^{\la,M})\vaa H^1(\gal,
 G_K^{\la}))
\]
can be identified with the set of
 $G_k^{\la}$-orbits of $\gal$-invariant sublattices
in the same $G_K^{\la}$-orbit.

Let $F$ be the stabilizer in $G$ of the space $W=\Omega M$.
$\Gamma$ its point-wise stabilizer. $H=F/\Gamma$.
There is a natural map
$G_K^{\la,M}\vaa H_K^M$,
which induces a map
\[
H^1(\gal,G_K^{\la,M})\vaa H^1(\gal,
 H_K^M)
\]
in cohomology.
If we are interested in lattices that are in the same $H_k$-orbit,
they will be classified by the kernel of the map
\[
H^1(\gal,G_K^{\la,M})\vaa H^1(\gal,
 H_K^M)\times H^1(\gal,G_K^{\la}).
\]
In the applications, $G$ is the stabilizer of a tensor $\tau$ of
type $(l,m)$, $W$ a subspace satisfying
$\tau(W^{\otimes l})\subseteq W^{\otimes m}$,
and $H=\stab_{GL(W)}(\tau|_W)$, where $\tau|_W$ is the restriction of
$\tau$ to $W$.
The condition on $W$ is vacuous if $m=0$.

\begin{ex}
Let $\tau=q$ is a quadratic form.
Inequivalent representations of $M_k$ by $\la_k$, that
become equivalent over $K$, are in correspondence with
\[
\ker\left(
H^1(\gal,\oink_{K,n}^{\la,M}(q))\vaa H^1(\gal,
\oink_{K,p}^M(q|_W))\times H^1(\gal,\oink_{K,n}^{\la}(q))
\right),
\]
where $n=\dim V$, $p=\dim W$,
and $q|_W$ is the restriction of $q$ to $W$.
\end{ex}

\begin{rmk}
All result in this paper apply also to lattices over rings of
$S$-integers. Absolute strong approximation must be replaced by strong
 approximation with respect to $S$.
\end{rmk}

\printindex

\end{document}